\documentclass{article}[12pt]
\usepackage{amsmath,amsfonts,amsthm}
\usepackage{amssymb}

\usepackage{epsfig}

\usepackage{graphics}
\usepackage{graphicx}

\usepackage{latexsym}
\usepackage{graphics}
\usepackage{amsmath}
\usepackage{amsthm}
\usepackage{xspace}
\usepackage{amssymb}
\usepackage{color}
\usepackage{cite}
\usepackage{latexsym}
\usepackage{graphics}
\usepackage{amsmath}
\usepackage{amsthm}
\usepackage{xspace}
\usepackage{amssymb}
\usepackage{epsfig}

\usepackage{enumitem}

\newcommand{\beql}[1]{\begin{equation}\label{#1}}
\newcommand{\eeql}{\end{equation}}
\newcommand{\eqn}[1]{(\ref{#1})}

\newcommand{\R}{\mathbb{R}}
\newcommand{\pr}{\mathbb{P}}
\newcommand{\E}{\mathbb{E}}

\def\P{{\mathbb P}}

\newcommand{\cx}{{\cal X}}

\newcommand{\cj}{{\cal J}}

\newcommand{\bI}{{\bf I}}

\newtheorem{thm}{Theorem}
\newtheorem{lem}[thm]{Lemma}

\newtheorem{cor}[thm]{Corollary}

\begin{document}

\title{Large-scale parallel server system with multi-component jobs 
%\\ (DRAFT)
}

\author
{
Seva Shneer \\
Heriot-Watt University\\
EH14 4AS Edinburgh, UK\\ 
\texttt{V.Shneer@hw.ac.uk}
\and
Alexander L. Stolyar \\
University of Illinois at Urbana-Champaign\\
%104 S. Mathews Avenue, Office 201C\\
Urbana, IL 61801, USA \\
\texttt{stolyar@illinois.edu}
}

\date{\today}

\maketitle

\begin{abstract}

A broad class of parallel server systems is considered, for which
we prove the steady-state asymptotic independence of server workloads, as the number of servers goes to infinity, while the system load remains sub-critical. Arriving jobs consist of multiple components. There are multiple job classes, and each class may be of one of two types, which determines the rule according to which the job components add workloads to the servers.
The model is broad enough to include as special cases some popular queueing models with redundancy, such as cancel-on-start and cancel-on-completion redundancy. 

Our analysis uses mean-field process representation and the corresponding mean-field limits.
In essence, our approach relies almost exclusively on three fundamental properties of the model: (a) monotonicity, (b) work conservation, (c) the property that, on average, ``new arriving workload prefers to go to servers with lower workloads.''

\end{abstract}


\iffalse
\noindent
{\em Abbreviated Title:} 
\fi


%\newpage

\section{Introduction}

In this paper we consider a broad class of parallel server systems, for which we prove the steady-state asymptotic independence of server workloads, as the number of servers goes to infinity, while the system load remains sub-critical. Our model is such that arriving jobs consist of multiple components. There are multiple job classes, and each class may be of one of two types. A job type determines the rule according to which the job components add workloads to the servers.
The model is broad enough to include as special cases some popular queueing models with redundancy, such as cancel-on-start and cancel-on-completion redundancy. 

More specifically, we consider a service system consisting of $n$ identical servers, processing work at unit rate. Jobs of multiple {\em classes} $j$ arrive according to independent Poisson processes of rate $\lambda_j n$. A job of each class $j$ consists of $k_j \ge 1$ {\em components}, while the $k_j$-dimensional exchangeable distribution $F_j$ determines the random component sizes, or workloads, $(\xi_1,\ldots, \xi_{k_j})$. 
(I.i.d. component sizes is a special case of exchangeability.)
Each class-$j$ job uniformly at random selects a subset of $d_j$ servers, $d_j \ge k_j$. Each job class may be of one of the two {\em types}, either {\em water-filling} or {\em least-load}. A job type determines the way in which the arriving job adds workload to the servers. For the least-load type
the component (random) workloads $(\xi_1,\ldots, \xi_{k_j})$ are added to the $k_j$ least-loaded servers out of the $d_j$ selected by the job.
For the water-filling type we describe the workload placement algorithm via the following illustration.
Suppose, $d_j=4$, $k_j=2$, the component sizes realization is $(10,5)$, and the workloads of the selected 4 servers are $5, 12, 7, 16$. Then, adding 10 units of the first component workload in the water-filling fashion brings the selected servers' workloads to $11, 12, 11, 16$. 
Before we place the next -- second -- component's workload, we exclude one of the selected servers that currently have the smallest workload -- it will be one of the servers with workload 11 in this illustration. 
Then, placing the second component's workload 5 in the water-filling fashion on the remaining 3
selected servers, brings the servers' workloads to $11, 14, 14, 16$. In general, after each component is placed,
the set of selected servers is reduced by excluding one of the servers with the smallest workload.

We assume that the system is sub-critically loaded, $\sum_j \lambda_j s_j < 1$, where $s_j$ is the total expected workload brought
by a class $j$ job. It is not hard to see that the system is stable for each $n$. Our {\bf main results,} Theorem~\ref{thm-main} and Corollary~\ref{cor-main}, prove the {\bf steady-state asymptotic independence} property: for any fixed integer $m\ge 1$, as $n\to\infty$, the steady-state workloads of a fixed set of $m$ servers (say, servers $1, \ldots,m$), become asymptotically independent. This property, in addition to be important in itself, in many cases allows one to obtain asymptotically exact system performance metrics, such as steady-state job delay mean or distribution. 

Our model is related to -- but not limited to -- queueing models with cancel-on-start and cancel-on-completion redundancy
\cite{Vulimiri-2013,shah2015redundant,gardner2015reducing,Harcol-Balter-2017,adan2018fcfs,ayesta2018unifying,ayesta2019redundancy}.
 In the redundancy models each job places its ``replicas'' on a selected subset of servers. The replicas may be served by their servers simultaneously. When a certain number of the job replicas start [resp., complete] their service, all other replicas are ``canceled'' and removed from the system. Hence the term cancel-on-start [resp., cancel-on-completion]. We postpone until the next section the detailed discussion of our model, including its relations to the models with redundancy. At this point we note that our least-load job type covers the cancel-on-start redundancy, and our water-filling job type covers the cancel-on-completion redundancy in the special case of i.i.d. exponentially distributed replica sizes. Moreover, our model allows multiple job classes, of different types.
 % within same system. 
 We also note that, for example, the model in \cite{Harcol-Balter-2017} and some of the models in \cite{hellemans2019performance}
 are special cases of ours; 
the steady-state asymptotic independence was used in those papers as a conjecture; our Corollary~\ref{cor-main} proves this conjecture (for those models in particular), thus formally substantiating the asymptotic steady-state performance metrics derived for those models in \cite{Harcol-Balter-2017}
and \cite{hellemans2019performance}.

Methodologically, this paper belongs to the line of work establishing the steady-state asymptotic independence in different contexts,
e.g. \cite{VDK96,BLP2012-jsq-asymp-indep, St2014_pull, St2015_pull}. Our approach is based 
on analyzing the mean-field (fluid) scaled process and its limit. One part of our analysis, namely establishing asymptotic independence of server workloads {\em over a finite time interval}, 
closely follows the previous work, namely \cite{GMP97,BLP2012-jsq-asymp-indep}. But, the main part of the analysis, namely the transition from the finite-interval asymptotic independence to the
steady-state asymptotic independence, is different from that in \cite{BLP2012-jsq-asymp-indep}. (Paper \cite{GMP97} does not have a steady-state asymptotic independence result.) Specifically, we rely on the dynamics -- transient behavior -- of the mean-field scaled process and its limit; 
in this sense, our approach is close to that in \cite{St2014_pull, St2015_pull}.
(The approach of \cite{BLP2012-jsq-asymp-indep} relies in an essential way on the direct steady-state estimates of the marginal workload distributions,
obtained in \cite{Bramson2011-jsq-stabil}.) At a high level, one may say that our approach relies almost exclusively on three fundamental properties of the model: (a) monotonicity, (b) work conservation, (c) the property that, on average, ``new arriving workload prefers to go to servers with lower workload.'' Consequently, we believe this approach for establishing the steady-state asymptotic independence may apply more broadly, to other models as well, as long as they possess these three fundamental properties.

The rest of the paper is organized as follows. A more detailed discussion of our model and results is given in Section~\ref{sec-discussion}, which is followed by a brief review of previous work in Section~\ref{sec-prev-work}. Section~\ref{sec-basic-notation} gives basic notation used throughout the paper. Section~\ref{sec-formal-results} presents our formal model and the main results, Theorem~\ref{thm-main} and Corollary~\ref{cor-main}. 
In Section~\ref{sec-more-general} we define some generalizations of our model and give their basic properties;
these generalizations, while may be of independent interest, are primarily for the purposes of the analysis of our original model.
Section~\ref{sec-auxilliary} contains more auxiliary facts used in the analysis.
Section~\ref{sec-asymp-indep-finite-int} gives results on the {\em finite-interval} asymptotic independence of the server workloads.
In Section~\ref{sec-fsp} we define limits of the mean-field (fluid) scaled processes; we call these limits fluid sample paths (FSP).
In Section~\ref{sec-fsp-from-empty} we study properties of the FSPs, starting specifically from ``empty'' initial state.
In Section~\ref{sec-fp} we define and study the FSP fixed point, which is the point to which the FSP trajectory, starting from the empty state, converges.
Finally, Section~\ref{sec-proof-main} contains the proof of the main result (Theorem~\ref{thm-main}), which employs the results developed
in Sections~\ref{sec-more-general}-\ref{sec-fp}.

\section{Discussion of the model and main results}
\label{sec-discussion}

The least-load job type is motivated by two scenarios. First, if we consider a system such that the current server workloads can be observed on a subset of $d_j$ selected servers, and a the job consists of $k_j$ components, it directly makes sense to place those components for service on the least-loaded $k_j$ of those $d_j$ servers. (See e.g. LL(d) policy in \cite{BLP2012-jsq-asymp-indep}, which is a special case of our model with a single least-load class with $d_j=d$, $k_j=1$.)
The second scenario arises in systems where the current workloads are not observable, and which use redundancy 
to improve performance. (See e.g. \cite{shah2015redundant} for a general motivation for redundancy.)
For example, suppose a class-$j$ job places $d_j$ ``replicas'' on $d_j$ randomly selected servers, and each server processes its work (replicas of different jobs) in the First-Come-First-Serve (FCFS) order. Suppose the job, to be completed, requires only $k_j$ (out of $d_j$) replicas to be processed, and as soon as the first $k_j$ replicas of the job {\em start} being processed, the remaining $d_j -k_j$ replicas of the job are ``canceled'' and immediately removed from their corresponding servers. This is usually referred to as {\em cancel-on-start} redundancy. (See e.g. \cite{ayesta2019redundancy,ayesta2018unifying}, where a special case $k_j=1$ is considered.)
We will call it {\em $(d_j,k_j)$-c.o.s. redundancy}, where $(d_j,k_j)$ are the parameters of class $j$.
Clearly, from the point of view of the server workload evolution (which need not be observable in this case), the described $(d_j,k_j)$-c.o.s. redundancy is equivalent to simply placing $k_j$ job replicas on the $k_j$ least loaded (out of $d_j$) servers, and not placing any workload on the remaining $d_j -k_j$ servers. Thus, a job class using $(d_j,k_j)$-c.o.s. redundancy can be equivalently viewed as a least-load job class in terms of our model,
with the $k_j$ job components being the first $k_j$ replicas. 

The water-filling job type motivation is also two-fold. First, suppose a job class $j$ is of the water-filling type, with $k_j=1$.
Consider a  class $j$ job;
it has one component. Suppose further that this component's workload can be arbitrarily divided between servers, in the sense that a parallel processing by multiple servers is allowed. (For example, the servers may represent different data transmission channels, with a job being a file that needs to be transmitted, and the job size being the file size.) Suppose the job can use $d_j$ randomly selected servers. The servers process workload in the FCFS order. The job is completed when all its workload is processed. Then, if the objective is to minimize the job completion time, its workload should be placed into the selected $d_j$ servers in the water-filling fashion. This can be done directly, if the workloads of the selected $d_j$ servers and the job workload are observable, or indirectly, as follows. The job joins the FCFS queues at each of the selected servers. When this job, at any of the selected servers, reaches the head of the queue -- i.e., can start using that server -- that server starts processing the job, possibly in parallel with other selected servers. The job is completed when the total amount of service it receives from all servers is equal to its size, at which point the job is removed from all queues. From the point of view of the server workload evolution (which need not be observable in this case), the described procedure is equivalent to simply placing the job's single component on $d_j$ selected servers in the water-filling fashion. 

The second motivation for the water-filling job type arises from {\em cancel-on-completion} redundancy (\cite{shah2015redundant,Vulimiri-2013,gardner2015reducing,Harcol-Balter-2017,ayesta2018unifying,adan2018fcfs,hellemans2019performance}). Suppose a class-$j$ job places $d_j$ job replicas on $d_j$ randomly selected servers. 
Each server processes its work (replicas of different jobs) in the FCFS order. Suppose the job, to be completed, requires only $k_j$ (out of $d_j$) replicas to be processed, and as soon as the first $k_j$ replicas of the job {\em complete} their service, the remaining $d_j -k_j$ replicas of the job are ``canceled'' and immediately removed from their corresponding servers. (Hence the name cancel-on-completion.) 
We will call this {\em $(d_j,k_j)$-c.o.c. redundancy}, where $(d_j,k_j)$ are the parameters of class $j$.
Suppose, in addition, that the {\em replica sizes for a class-$j$ job are i.i.d. exponential} random variables with mean $s_j/k_j$.
(This additional assumption, as well as the further assumption that $k_j=1$, is used, e.g., in \cite{Harcol-Balter-2017}.) 
Under this additional assumption (of i.i.d. exponential replica sizes), 
it is easy to observe that, from the point of view of the servers' workload evolution (which need not be observable), the described $(d_j,k_j)$-c.o.c. redundancy is equivalent to placing on the selected $d_j$ servers a water-filling job with the following parameters: $(d_j,k_j)$ are the same as above, and the component sizes are i.i.d. exponential random variables with mean $s_j/k_j$. Indeed, the job component 1 places (stochastically) exactly the same amounts of additional workload on the servers as the workloads placed by all replicas {\em up to the time of the first replica service completion.} Similarly, the job component 2 places (stochastically) exactly the same amounts of additional workload on the servers as the workloads placed by all replicas {\em from the time of the first replica service completion until the time of the second replica service completion}. And so on. Thus, a job class using $(d_j,k_j)$-c.o.c. redundancy (under the additional assumption of i.i.d. exponential replica sizes)
can be equivalently viewed as a water-filling job class in terms of our model,
with parameters $(d_j,k_j)$  and i.i.d. exponentially distributed components.

We see that our model is very broad. In \cite{BLP2012-jsq-asymp-indep} it is proved
(among other results) that the steady-state asymptotic 
independence (our Corollary~\ref{cor-main}) holds for the special case of a single, least-load job class with $k_j=1$.
(See LL(d) model in \cite{BLP2012-jsq-asymp-indep}.) The $(d,1)$-c.o.c. model with  i.i.d. exponential replica sizes,
considered in \cite{Harcol-Balter-2017}, is a special case of our model, with a single, water-filling job class, with $k_j=1$ and i.i.d. exponential component sizes. One of the models considered in \cite{hellemans2019performance} (called LL(d,k,0) there) 
is a special case of ours, with a single least-load job class. In both \cite{Harcol-Balter-2017} and \cite{hellemans2019performance}
the steady-state asymptotic independence was used as a conjecture; our Corollary~\ref{cor-main} proves this conjecture for both models.
Furthermore, since our model allows multiple job classes of different types, Corollary~\ref{cor-main} establishes 
the steady-state asymptotic independence, for example, for a system with two job classes -- one as in \cite{Harcol-Balter-2017} and
one as the LL(d,k,0) class in \cite{hellemans2019performance}.

\section{Previous work}
\label{sec-prev-work}

The work on the steady-state asymptotic independence in the large-scale regime, with the number of servers and the arrival rate increasing to infinity, while the system load remains sub-critical, includes, e.g., papers \cite{VDK96,BLP2012-jsq-asymp-indep, St2014_pull}.
Papers \cite{VDK96,BLP2012-jsq-asymp-indep} prove this for the celebrated ``power-of-d'' choices algorithm, where each arriving (single-component) job joins the shortest queue out of $d$ randomly selected; \cite{VDK96} does this for the exponentially distributed job sizes, while 
 \cite{BLP2012-jsq-asymp-indep} extends the results to more general job size distributions, namely those with decreasing hazard rate (DHR). 
 Note that a standard power-of-d choices algorithm is {\em not} within the framework of our model, because job placement decisions depend on the queue lengths (numbers of jobs), as opposed to depending on the workloads. However, \cite{BLP2012-jsq-asymp-indep} also considers -- and establishes the steady-state asymptotic independence for -- the LL(d) model, which is a special case of our model with the single, least-load job class with $d_j=d$ and $k_j=1$. Note that, equivalently, this is the single-class $(d,1)$-c.o.s. redundancy model.
 The main results of \cite{BLP2012-jsq-asymp-indep} in turn rely on the uniform estimates of the marginal stationary distribution of a single server state, obtained in \cite{Bramson2011-jsq-stabil}. Paper \cite{St2014_pull} proves the steady-state asymptotic independence under a pull-based algorithm, also for the model with single-component jobs, having DHR size distributions. (The model in \cite{St2014_pull} is also {\em not} within the framework of the present paper model.)
 
 For the redundancy models, such as in \cite{Vulimiri-2013,shah2015redundant,gardner2015reducing,Harcol-Balter-2017,adan2018fcfs,ayesta2018unifying,ayesta2019redundancy,hellemans2019performance}, we are not aware of any prior steady-state asymptotic independence
 results, besides the already mentioned $(d,1)$-c.o.s. result in \cite{BLP2012-jsq-asymp-indep}. However, the steady-state asymptotic independence {\em conjecture} is often used (e.g. \cite{Vulimiri-2013,Harcol-Balter-2017,hellemans2019performance}) 
 to obtain estimates of the steady-state performance metrics of large scale systems.
 
 Paper \cite{Vulimiri-2013} introduces redundancy as a way to reduce job delays. It considers the $(d,1)$-c.o.c. redundancy model, with generally distributed replica sizes. (As such, this model is {\em not} within the framework of our model.) The paper uses the steady-state asymptotic independence {\em conjecture} to estimate the average job delay when the system is large.
 
 Paper \cite{shah2015redundant} introduces and motivates the $(d,k)$-c.o.c. redundancy model, and establishes a variety of monotonicity properties of the average job delay with respect to the selection set size $d$, under different assumptions on the replica size distribution. 
 Some of the results of \cite{shah2015redundant} are for the $(d,k)$-c.o.c. redundancy model with i.i.d. exponential replica sizes,
 which is a special case of our model, but \cite{shah2015redundant} does not consider the asymptotic regime with $n\to\infty$.

As already described earlier, paper \cite{Harcol-Balter-2017} studies the $(d,1)$-c.o.c. redundancy model with i.i.d. exponential replica sizes,
and obtains asymptotically exact expressions for the job delay distribution, based on the steady-state asymptotic independence {\em conjecture}. Our results prove this conjecture, thus completing formal substantiation of those asymptotic expressions.

Paper \cite{hellemans2019performance} studies general redundancy models -- more general than c.o.s. and c.o.c. that we described earlier -- and uses the steady-state asymptotic independence {\em conjecture} to characterize and compute steady-state performance metrics. Some 
(not all) of the redundancy schemes in \cite{hellemans2019performance} are within our model framework. For example, LL(d,k,0) redundancy in 
\cite{hellemans2019performance} is a special case of our model with a single least-load class. Thus, again, 
by proving the steady-state asymptotic independence,
our results 
complete formal substantiation of some of the asymptotic results in \cite{hellemans2019performance}.

Papers \cite{adan2018fcfs,ayesta2018unifying,ayesta2019redundancy} derive explicit product-form stationary distributions for 
different versions of $(d,1)$-c.o.c. and $(d,1)$-c.o.s. redundancy, assuming i.i.d. exponential replica sizes.

\section{Basic notation}
\label{sec-basic-notation}

We denote by $\R$ and $\R_+$ the sets of real and real non-negative numbers, respectively, and by $\R^n$ and $\R_+^n$ the corresponding $n$-dimensional product sets. By $\bar \R_+ \doteq \R_+ \cup \{\infty\}$ we denote the one-point compactification of $\R_+$, where 
$\infty$ is the point at infinity, with the  natural topology.
We say that a function is RCLL if it is {\em right-continuous with left-limits}. Inequalities applied to vectors [resp.  functions] are understood component-wise [resp. for every value of the argument]. The sup-norm of a scalar function $f(w)$ is denoted $\|f(\cdot)\| \doteq \sup_w |f(w)|$; the corresponding convergence is denoted by $\stackrel{u}{\rightarrow}$. {\em U.o.c.} convergence means {\em uniform on compact sets} convergence, and is denoted by $\stackrel{u.o.c.}{\rightarrow}$. We use notation: $a\vee b \doteq \max\{a,b\}$, 
$a\wedge b \doteq \min\{a,b\}$.
Abbreviation WLOG means {\em without loss of generality}. 

For a random process $Y(t), ~t\ge 0,$ we denote by $Y(\infty)$ the random value of $Y(t)$ in a stationary regime (which will be clear from the context). Symbol $\Rightarrow$ signifies convergence of random elements in distribution; $\stackrel{P}{\rightarrow}$ means convergence in probability.
 {\em W.p.1} means {\em with probability one.}
{\em I.i.d.} means {\em independent identically distributed.}
Indicator of event or condition $B$ is denoted by $\bI(B)$. If $X,Y$ are random elements taking values in set $\cx$, on which a partial order $\le$ is defined, then the stochastic order
$X \le_{st} Y$ means that $X$ and $Y$ can be coupled (constructed on a common probability space) so that $X \le Y$ w.p.1.

We will use the following non-standard notation. Suppose $f^n_w$, $n\to\infty$, is a sequence of random functions of $w$, and $f_w$ is a 
deterministic function of $w$. Then, for a fixed $w$,
\beql{eq-P-liminf}
(P) \!\liminf_{n\to\infty} f^n_w \ge f_w ~~~\mbox{means} ~~~ [(f^n_w - f_w) \wedge 0] \stackrel{P}{\rightarrow} 0, ~n\to\infty,
\end{equation}
and for a subset $A$ of the domain of $w$,
\beql{eq-P-liminf2}
(P) \!\liminf_{n\to\infty} (f^n_w, ~w\in A) \ge (f_w, ~w\in A) ~~~\mbox{means} ~~~ \inf_{w\in A} [(f^n_w - f_w) \wedge 0] \stackrel{P}{\rightarrow} 0, ~n\to\infty.
\end{equation}
Analogously,
$$
(P) \!\limsup_{n\to\infty} f^n_w \le f_w ~~~\mbox{means} ~~~ [(f^n_w - f_w) \vee 0] \stackrel{P}{\rightarrow} 0, ~n\to\infty,
$$
$$
(P) \!\limsup_{n\to\infty} (f^n_w, ~w\in A) \le (f_w, ~w\in A) ~~~\mbox{means} ~~~ \sup_{w\in A} [(f^n_w - f_w) \vee 0] \stackrel{P}{\rightarrow} 0, ~n\to\infty,
$$
$$
(P) \!\lim_{n\to\infty} f^n_w = f_w ~~~\mbox{means} ~~~ f^n_w  \stackrel{P}{\rightarrow} f_w, ~n\to\infty,
$$
$$
(P) \!\lim_{n\to\infty} (f^n_w, ~w\in A) = (f_w, ~w\in A) ~~~\mbox{means} ~~~ \sup_{w\in A} |f^n_w - f_w| \stackrel{P}{\rightarrow} 0, ~n\to\infty.
$$

\section{Formal model and main results}
\label{sec-formal-results}

\subsection{Model}

There are $n$ identical servers.  The unfinished work of a server at a given time will be referred to as its {\em workload.}
Each server processes its workload at rate $1$. 
There is a finite set $\cj$
of job classes, which are indexed by $j$. (Set $\cj$ does not depend on $n$.) Jobs of class $j$ arrive according to a Poisson process of rate $\lambda_j(n) n$.
Associated with each class $j$ there are three parameters: integers $k_j$ and $d_j$
such that $1\le k_j \le d_j$, and the exchangeable probability distribution $F_j$ on $\R_+^{k_j}$. 
A class-$j$ job consists of $k_j$ {\em components}, with each component having a (random) {\em size} (which is the amount of new workload 
this component brings); $F_j$ is the joint distribution of random component sizes $(\xi_1^{(j)},\ldots, \xi_{k_j}^{(j)})$ for a class-$j$ job. 
Exchangeability of $F_j$ means that it is invariant w.r.t. permutations of component indices.
We assume that 
$s_j \doteq \E \sum_\ell \xi_\ell^{(j)} = k_j \E \xi_1^{(j)} < \infty$.
WLOG, we can and do assume that $s_j >0$ for each class $j$.
We will denote $d\doteq \max_j d_j$.

Each job class $j$ may be of one of the two {\em types}, either {\em water-filling} or {\em least-load}. 
The corresponding non-intersecting subsets of $\cj$ we denote by $\cj_{wf}$ and $\cj_{ll}$. (Either of them may be empty.)
A job type determines the way in which the arriving job adds workload to the servers. We will describe the job types separately.

{\em A least-load job class $j\in \cj_{ll}$.} When such a job arrives, $d_j$ servers are selected uniformly at random; these servers form the {\em selection set} of the job. Then $k_j$ of the selected servers, that are least-loaded (have the smallest workload), are picked;
the workload ties are broken in an arbitrary fashion. Then, independently of the process history, random component sizes $(\xi_1,\ldots, \xi_{k_j})$ are drawn according to distribution $F_j$. 
Then, workload $\xi_1$ is added to the least-loaded of those servers, $\xi_2$ is added to the second least-loaded of those servers, and so on.

{\em A water-filling job class $j\in \cj_{wf}$.} When such a job arrives, its {\em selection set} of $d_j$ servers is selected uniformly at random. Then, independently of the process history, random component sizes $(\xi_1,\ldots, \xi_{k_j})$ are drawn according to distribution $F_j$. 
We ``take'' the first component, and place its $\xi_1$-size workload on the servers within the selection set in the ``water-filling'' fashion. 
(For example, suppose the selection set consists of $4$ servers, $1,2,3,4$, with workloads $W_1=5, W_2=12, W_3=7, W_4=16$,
and suppose $\xi_1=10$. Then, adding the workload of size $10$ to these servers in the water-filling fashion will result in the new workloads
being $W_1=11, W_2=12, W_3=11, W_4=16$. That is servers 1 and 3 will receive non-zero additional workloads, $6$ and $4$, respectively, and
will end up with equal workload $11$. Servers  2 and 4 will not receive any of the first component's workload.)
After this, there will be one or more selected servers that currently have the smallest workload.
(Servers 1 and 3 in the illustration above.)
Let us call them component-1 servers.
Then we pick one of the component-1 servers (in an arbitrary fashion), and {\em exclude it} from further workload placement by this job. Then, 
we ``take'' the second component, and place its $\xi_2$-size workload on the remaining $d_j-1$ servers by continuing the water-filling. 
The servers that will have the smallest workload after that we call component-2 servers.
Then we exclude one of the component-2 servers, and so on, until the workload of all  $k_j$ components is placed.
(Note that we could define an additional -- different -- water-filling type, such that the water-filling continues to use all $d_j$ selected servers, without excluding one of the servers after each component placement.
This, however, is just a special case of the type we just defined, with $k_j$ components replaced by the single component of the
size $\sum_i \xi_i$.)

By the model definition, for each class $j$, regardless of its type, the total expected additional workload it brings to the system is
equal to $s_j$.

 \subsection{Asymptotic regime. Mean-field scaled process} 

 We consider the sequence of systems with $n\to\infty$, and assume 
 $$
 \lambda_j(n) \to \lambda_j >0, ~~j \in \cj.
 $$
 Further assume that the system is (asymptotically) sub-critically loaded
\begin{equation} 
\label{eq:def_rho}
\rho \doteq \sum_j \lambda_j s_j < 1.
\end{equation}

Denote the (limiting) total job arrival rate per server by
\begin{equation} 
\label{eq:def_lambda}
\lambda \doteq \sum_j \lambda_j.
\end{equation}
WLOG, we can and will assume $\lambda=\rho$. (We can achieve this by rescaling time, if necessary.)

To improve paper readability, let us assume that $\lambda_j(n) = \lambda_j$ for each $n$. Having converging $\lambda_j(n)$ does not change anything of substance, but clogs exposition. (However, we do need and will use the fact that our results hold for converging arrival rates.) Similarly, throughout the paper, we will often consider ``$an$ servers'' for some real $a$, ignoring the fact that $an$ may be non-integer; it would be more precise to consider, for example,  ``$an$-rounded-up servers,'' but it would just clog the exposition, rather than creating any difficulties. 

From now on, the upper index $n$ of a variable/quantity will indicate that it pertains to the system with $n$ servers,
or $n$-th system.
Let $W_i^n(t)$ denote the workload
of server $i$ at time $t$ in the $n$-th system. (When $W_i^n(t)=0$ we say that server $i$ at time $t$ is empty.)
Consider the following {\em mean-field}, or {\em fluid}, scaled quantities:
\beql{eq-x-def2}
x^n_w(t) \doteq (1/n) \sum_i \bI\{W_i^n(t)> w\}, ~~ w \ge 0.
\end{equation}
That is, $x^n_w(t)$ is the fraction of servers $i$ with $W_i^n(t)> w$.
Then
 $x^n(t)=(x^n_w(t), ~w\ge 0)$
is the system state at time $t$;
note that $x^n_0(t)$ 
%$\rho^n(t) \equiv  x^n_0(t)$ 
is the fraction of busy servers (the instantaneous system load).

For any $n$, the state space of the process $(x^n(t), ~t\ge 0)$ is a subset of a common (for all $n$) 
state space $\cx$, whose elements
$x=(x_w, ~w\ge 0)$ are non-increasing RCLL functions of $w$, 
with values $x_w \in [0,1]$. 
An element $x \in \cx$ defines a probability measure on $\bar \R_+$, with $1-x_w$ being the measure of $[0,w]$ for $0\le w < \infty$.
Denote $x_\infty \doteq \lim_{w\to\infty} x_w$; then $x_\infty$ is the measure of $\{\infty\}$.
An element $x \in \cx$ we will call {\em proper}, if $x_\infty = 0$, i.e. if the corresponding probability measure is concentrated 
on $\R_+$. We will equip the space $\cx$ with the topology of weak convergence of measures on $\bar \R_+$;
equivalently, $y \to x$ if and only if $y_w \to x_w$ for each $0 < w < \infty$ where $x$ is continuous.
We also can and do equip $\cx$ with a metric consistent with the topology. 
Obviously, $\cx$ is compact.

For any $n$, process $x^n(t), ~t\ge 0,$ is Markov with state space $\cx$,
and with sample paths being RCLL functions of $t\ge 0$.
% with values in $\cx$. 
Moreover, this is a renewal process, with renewals occurring when all servers become empty. 

Under the subcriticality assumption \eqn{eq:def_rho}, i.e. $\rho=\lambda<1$, 
the stability (positive Harris recurrence) of the process $(x^n(t), ~t\ge 0)$, for any $n$, is not hard to establish.
(Positive recurrence in this case simply means 
that the expected time to return to the empty state is finite.)
It can be established, for example, using the fluid limit technique, analogously to the way it is done in \cite{Foss-Chernova-1998}. 
The key property that the fluid limit for our model shares with that in \cite{Foss-Chernova-1998} is that if there is a subset of servers, whose fluid workloads are greater than in the rest of the servers, the average per-server rate at which the servers within the subset will receive new workload is at most $\rho<1$. (See \eqn{eq-rho_A-prelim} in Section~\ref{sec-equiv}.) We do not provide further details of the stability proof.

Given that the process $(x^n(t), ~t\ge 0)$ is stable, it has a unique stationary distribution.
Let $x^n(\infty)$ be a random element whose distribution is the stationary distribution of the process; in other words,
this is a random system state in the stationary regime.

\subsection{Main results}

\begin{thm}
\label{thm-main}
There exists a 
unique 
proper element $x^* \in \cx$, with $x^*_0 = \lambda = \rho$, such that 
\beql{eq-conv-main}
x^n(\infty) \Rightarrow x^*, ~~n\to\infty.
\eeql
Function $x^*_w, ~w\ge 0,$ is Lipschitz continuous and strictly decreasing (and then everywhere positive).
\end{thm}

\begin{cor} [Steady-state asymptotic independence]
\label{cor-main}
For any fixed integer $m \ge 1$, the following holds. For each $n$, 
let $[W^n_1(\infty), \ldots, W^n_m(\infty)]$ denote the random value of $[W^n_1(t), \ldots, W^n_m(t)]$ in the stationary regime.
Then
%, as $n\to\infty$,
\beql{eq-indep-main}
[W^n_1(\infty), \ldots, W^n_m(\infty)] \Rightarrow [W^*_1, \ldots, W^*_m], ~~ n\to\infty,
\eeql
where random variables $W^*_i$ are i.i.d., with
$\P\{W^*_i > w\} = x^*_w, ~w\ge 0$.
\end{cor}

Corollary~\ref{cor-main} follows from Theorem~\ref{thm-main} and the symmetry between (exchangeability of) servers.
Indeed, by the symmetry, the distribution of $[W^n_1(\infty), \ldots, W^n_m(\infty)]$ is equal to the joint distribution of workloads of $m$ servers
chosen uniformly at random. Since $m$ is fixed and as $n\to\infty$ the steady-state empirical measure $x^n(\infty)$ of 
the server workloads converges to the deterministic element $x^*$, the statement easily follows. 
%It is easy to see that, in fact, \eqn{eq-indep-main} is equivalent to \eqn{eq-conv-main}.

\section{More general systems}
\label{sec-more-general}

\subsection{Infinite workloads and truncation. Monotonicity properties.}
\label{sec-monotonicity-etc}

For the purposes of our analysis, it will be convenient to consider two generalizations of our model. (These more general systems may be of independent interest as well.)

First, we generalize our original system defined above, by allowing some of the servers to have an infinite workload. Specifically, if server $i$ workload is initially infinite,
$W_i^n(0)=\infty$, then, by convention, it remains infinite at all times, $W_i^n(t)=\infty$, $t\ge 0$. The same workload placement rules apply even if some server workloads are infinite, with the convention that an infinite workload remains infinite when ``more'' workload is added to it.
Note that if one or more server workloads are initially infinite, this implies that $x^n_\infty(0)>0$ and $x^n_\infty(t)=x^n_\infty(0)$ for all $t\ge 0$.

A second convenient generalization is a system, where the workload of the servers is truncated at some level $c$, where $0 \le c \le \infty$. Such a truncated system is defined exactly as the original one, except when an arriving job adds workload to servers, each server's workload is capped (truncated) at level $c$ every time the algorithm would increase it above $c$. The workload lost due to truncation is removed from the system. The case $c=\infty$ corresponds to the original, non-truncated system, where the arriving workload is never lost. Note that, if $c<\infty$, then the stability for any $n$ (and any $\lambda$) is automatic.
The process corresponding to the truncated system with parameter $c$, we denote by $x^{n,c}(t), ~t\ge 0;$ if superscript $c$ is absent, 
this corresponds to $c=\infty$, i.e. the process $x^{n}(\cdot)$ is for the original non-truncated system.

Finally, if the process starts specifically from the ``empty'' initial state (with all servers having zero initial workload), we will add superscript $\emptyset$ to the process notation: 
$x^{n,c,\emptyset}(t), ~t\ge 0$; therefore, $x^{n,c,\emptyset}_0(0)=0$. So, for example, $x^{n,\emptyset}(\cdot)$
denotes the original non-truncated process, starting from the empty state.

The analysis in this paper relies on the system monotonicity, and related properties. We will need several such properties. They are all related and rather simple. 

\begin{lem}
\label{lem-monotone1}
Consider two versions of the process, $x^{n,c}(\cdot)$ and $\hat x^{n, \hat c}(\cdot)$,
such that $x^{n,c}(0) \le \hat x^{n, \hat c}(0)$, $c \le \hat c$.
Then these processes can be coupled so that, w.p.1,
\beql{eq-key-monotone}
x^{n,c}(t) \le \hat x^{n, \hat c}(t), ~\forall t\ge 0.
\eeql
Furthermore, if the process $\hat x^{n, \hat c}(\cdot)$ is modified so that, in addition to the job arrival process (as defined in our model), arbitrary amounts of workload may be added at arbitrary times to arbitrary servers, the property \eqn{eq-key-monotone} still holds.
\end{lem}

 {\em Proof.}  As far as the mean-field scaled processes $x^{n,c}(\cdot)$ and $\hat x^{n, \hat c}(\cdot)$
 are concerned, WLOG, we can assume that, after each job arrival and/or other workload addition(s), the actual servers $1,\ldots,n$ are relabeled, so that the workloads
 $W_1^n, \ldots, W_n^n$ are non-decreasing. Then, for the two processes it is sufficient to couple in the natural way
 the arrival processes and the job selection sets, to see that \eqn{eq-key-monotone} must prevail at all times.
$\Box$

From Lemma~\ref{lem-monotone1}, we obtain the following

\begin{cor} 
\label{cor:zero}
For any $0 \le c \le \infty$, the process $x^{n,c,\emptyset}(\cdot)$ is monotone in time $t \ge 0$, namely
$$
x^{n,c,\emptyset}(t_1) \le_{st} x^{n,c,\emptyset}(t_2), ~~\forall t_1 \le t_2.
$$
\end{cor}

\begin{lem}
\label{lem-equality-persists}
Consider two versions of the process, $x^{n,c}(\cdot)$ and $\hat x^{n, \hat c}(\cdot)$,
such that $c, \hat c \in [0, \infty]$. Suppose that for some fixed $w \in [0,c\wedge \hat c]$,
we have $x_v^{n,c}(0) = \hat x_v^{n, \hat c}(0)$ for $0\le v \le w$.
Then these processes can be coupled so that, w.p.1, for $t\in [0,w]$ and $v\in [0,w-t]$,
\beql{eq-persists}
x_v^{n,c}(t) = \hat x_v^{n, \hat c}(t).
\eeql
\end{lem}

{\em Proof.}  We couple the two processes in the natural way, as in the proof of Lemma~\ref{lem-monotone1}. 
The proof then follows by induction on the times of job arrivals in the interval $[0,w]$. 
Indeed, if $t_1\le w$ in the time of the first job arrival, \eqn{eq-persists} of course holds for all $t\in [0,t_1)$.
It is then easy to see that the changes of $x_v^{n,c}$ and  $\hat x_v^{n, \hat c}$ for $v\le w-t_1$, at time $t_1$, 
only depend on those servers with workloads at most $w-t_1$, which are the same for both systems; we also observe that 
if any of those servers changes its workload to a value not exceeding $w-t_1$, the change will be exactly same 
in both systems. Then \eqn{eq-persists} holds for $t=t_1$. Then, \eqn{eq-persists} holds until the time $t_2$ of the second job arrival
or $w$, whichever is smaller. And so on.
$\Box$

Lemma~\ref{lem-equality-persists} and Lemma~\ref{lem-monotone1} imply the following more general form of Lemma~\ref{lem-equality-persists}.

\begin{lem}
\label{lem-equality-persists-gen}
Consider two versions of the process, $x^{n,c}(\cdot)$ and $\hat x^{n, \hat c}(\cdot)$,
such that $c, \hat c \in [0, \infty]$. Suppose that for some fixed $w \in [0,c\wedge \hat c]$,
we have $x_v^{n,c}(0) \le \hat x_v^{n, \hat c}(0)$ for $0\le v \le w$.
Then these processes can be coupled so that, w.p.1, for $t\in [0,w]$ and $v\in [0,w-t]$,
\beql{eq-persists-gen}
x_v^{n,c}(t) \le \hat x_v^{n, \hat c}(t).
\eeql
\end{lem}

\subsection{Equivalent representation of a system with some workloads being infinite.}
\label{sec-equiv}

Let $b\in (0,1]$ be fixed. For each $n$, consider the system with the initial state such that
$(1-b)n$ servers have infinite workloads, while the remaining $bn$ servers' workloads are finite. Let $B=B(n)$ denote the set of servers with finite workloads. Then, for each $n$, the evolution of the subsystem consisting 
of servers in $B$ -- let us call it $B$-subsystem -- can be equivalently described as follows. The number of servers is $bn$.
Each job class $j$ ``breaks down'' into multiple classes $(j,m)$, $m=1,\ldots,d_j$, as follows. Let $\pi_{j,m}(n)$ be the probability that exactly $m$ servers selected by a class $j$ job, will be in $B$. Note that
$$
\sum_m \pi_{j,m}(n) m = b d_j.
$$
Then, for a given $n$, class $(j,m)$ in the $B$-system has the following parameters: arrival rate per server $\lambda_{j,m}(n)
=\lambda_j \pi_{j,m}(n)/b$, $d_{j,m} = m$, $k_{j,m}= k_j \wedge m$, the distribution $F_{j,m}$ of the component sizes is 
the projection of the distribution $F_j$ on the first $m$ components. 
($d_{j,m}$, $k_{j,m}$ and $F_{j,m}$ do not depend on $n$.)
Clearly, as far as evolution of the $B$-system is concerned,
this new description is consistent with the actual behavior. The load of the $B$-system is
$$
\rho_B(n) = \sum_{j,m} \lambda_{j,m}(n) s_j (k_{j,m}/k_j). 
$$
Recall that the load of the original system, for any $n$, is $\rho = \sum_j \lambda_j s_j$.

The following fact is very intuitive -- by the nature of the workload placement algorithm, the arriving workload ``prefers'' servers with finite workloads.

\begin{lem}
\label{lem-B-load}
For each $n$,
\beql{eq-rho_B-prelim}
\rho_B(n) \ge \rho.
\eeql
\end{lem}

{\em Proof.}  We can write:
$$
\rho = \sum_j \lambda_j s_j = \sum_j \sum_m  \frac{\lambda_j \pi_{j,m}(n)}{b}  \frac{m}{d_j}  s_j
\le \sum_j \sum_m  \frac{\lambda_j \pi_{j,m}(n)}{b}  \frac{m \wedge k_j}{d_j \wedge k_j} s_j
= \sum_j \sum_m  \lambda_{j,m}(n)  \frac{k_{j,m}}{k_{j}}  s_j = \rho_B(n).
$$
$\Box$

Note that,  if $\rho_A(n)$ is the load of the complementary subsystem, consisting of the $(1-b)n$ infinite-workload servers,  then
$\rho = b \rho_B(n) + (1-b)\rho_A(n)$ and, therefore,
\beql{eq-rho_A-prelim}
\rho_A(n) \le \rho.
\eeql

Consider now a sequence of the above systems, with $n\to\infty$. Recall that the number of servers in the $B$-system is $bn$.
Note that 
$$
\lim_n \pi_{j,m}(n) =  \pi_{j,m} = \frac{d_j !}{m! (d_j-m)!} b^m (1-b)^{d_j-m}.
$$
Then,
$$
\lambda_{j,m}(n) \to \lambda_{j,m} = \lambda_j \pi_{j,m}/b,
$$
and the $B$-subsystem (limiting) load is
\beql{eq-rho_B}
\rho_B = \lim_n \rho_B(n) 
\ge \rho. 
\eeql
We see that the sequence of $B$-systems is just like our original sequence of systems, but has different parameters. 
(Recall that our original model does allow converging arrival rates per server, not just constant.)

\section{Some auxiliary facts}
\label{sec-auxilliary}

\begin{lem}
\label{lem-exp-decay-lower-bound}
Let $a\in [0,1]$ be fixed. 
Consider a sequence of processes such that, for each $n$, at time $0$, we identify a subset, consisting of $a n$ servers.
As the process evolves, for $t\ge 0$, we will keep track of those servers -- let us call them ``tagged.''
Denote by $f^n(t)$, $t\ge 0$, the (scaled) number of the tagged servers, which are not selected by any new job arrival in the interval $[0,t]$.
Then, for any fixed $t\ge 0$,
\beql{eq-exp-decay} 
(P) \!\liminf_{n\to\infty} f^n(t) \ge a e^{-\lambda d^2 t}.
\eeql
Since, by definition, $f^n(t)$ is non-increasing in $t$, as a corollary of \eqn{eq-exp-decay}, we 
obtain the following stronger property: for any $t\ge 0$,
\beql{eq-exp-decay-unif} 
(P) \!\liminf_{n\to\infty} (f^n(\tau), ~\tau \in [0,t]) \ge a e^{-\lambda d^2 t}.
\eeql
\end{lem}

 {\em Proof.}  Using coupling, we see that the stochastic lower bound $\hat f^n(t)$ of the process $f^n(t)$ can be obtained by considering the following ``worst case'' unaffected tagged set scenario: (a) {\em each} new job arrival selection set consists of  $d$ servers and (b) if at least one of the selected servers is within the set of currently unaffected tagged servers, the latter set is reduced by $d$ servers. 
 For the worst case unaffected tagged set, $\hat f(t)=a e^{-\lambda d^2 t}$ is the deterministic mean-field (fluid) limit,
 solving $\hat f'(t) = -(\lambda d) d \hat f(t)$ with $\hat f(0)=a$;
 here $\lambda d \hat f(t)$ is the (scaled, limiting) rate at which arriving jobs select a server within the set and $d$ is the number of servers removed upon each such event.
 Namely, using standard techniques
 (``large number of servers'' fluid limit), cf. \cite{PTW}, it is straightforward to show that,
 for any $t\ge 0$,
$$
(P) \!\lim_{n\to\infty} (\hat f^n(\tau), ~\tau \in [0,t]) = (\hat f(\tau), ~\tau \in [0,t]),
$$
which then implies \eqn{eq-exp-decay}.
 $\Box$
 
\begin{lem}
\label{lem-jump-remains}
Let $a\in [0,1]$ and $h>0$ be fixed. 
Consider a sequence of processes $x^n(\cdot)$ with initial states $x^n(0)$ satisfying the following condition: the (scaled) number
of servers with workload {\em exactly equal} to $h$, is at least $a$; namely, $x^n_{h-}(0) - x^n_{h}(0) \ge a$.
Then,
\beql{eq-jump-remains} 
(P) \!\liminf_{n\to\infty} (x^n_{h-t-}(t) - x^n_{h-t}(t), ~0 \le t < h) \ge a e^{-\lambda d^2 h}.
\eeql
(Informally, in words, ``when $n$ is large, then with high probability a positive, bounded away from zero, jump in $x^n(t)$ ``moves'' left at speed $1$ from initial point $w=h$.)
\end{lem}

 {\em Proof.}  Consider the servers with initial workload exactly equal to $h$ as tagged servers, and apply Lemma~\ref{lem-exp-decay-lower-bound}.
 $\Box$

 \section{Asymptotic independence over a finite interval}
 \label{sec-asymp-indep-finite-int}
 
 The constructions and the results in this section closely follow those in  \cite{BLP2012-jsq-asymp-indep} (Section 7) and (to a lesser degree) 
 in \cite{GMP97} (proofs of Lemmas 3.1 and 3.2). We give them here (along with  proofs) in the setting/notation that we need for our model.

Suppose a finite set of fractions (a probability distribution) $a_1, a_2, \ldots, a_K$ is fixed, where all $a_k>0$ and $\sum_k a_k=1$.
Also fixed is a set of numbers $w_k \in [0,\infty]$, $k=1,\ldots,K$. Let the truncation parameter $c\in [0,\infty]$ be fixed.
In this section, we consider a sequence of our systems, 
indexed by $n\to\infty$, with initial states such that $a_k n$ servers have workload exactly $w_k$, $k=1,\ldots,K$. Suppose, initially the server indices $1,\ldots,n$ are assigned in the order of a server set permutation chosen uniformly at random. This means, in particular, that $\pr\{W^n_1(0) = w_k\} = a_k$. 

We now formally construct a random process $U_1(t), ~t\ge 0$. Lemma~\ref{lem-graph-rep2} below will show that, for each $t$, $W_1^n(t) \Rightarrow U_1(t)$
as $n\to\infty$. So, informally speaking, this is a construction of the evolution of a server workload in a system
with ``infinite number of servers.'' 

\begin{figure}
\centering
\includegraphics[width=1.0\linewidth]{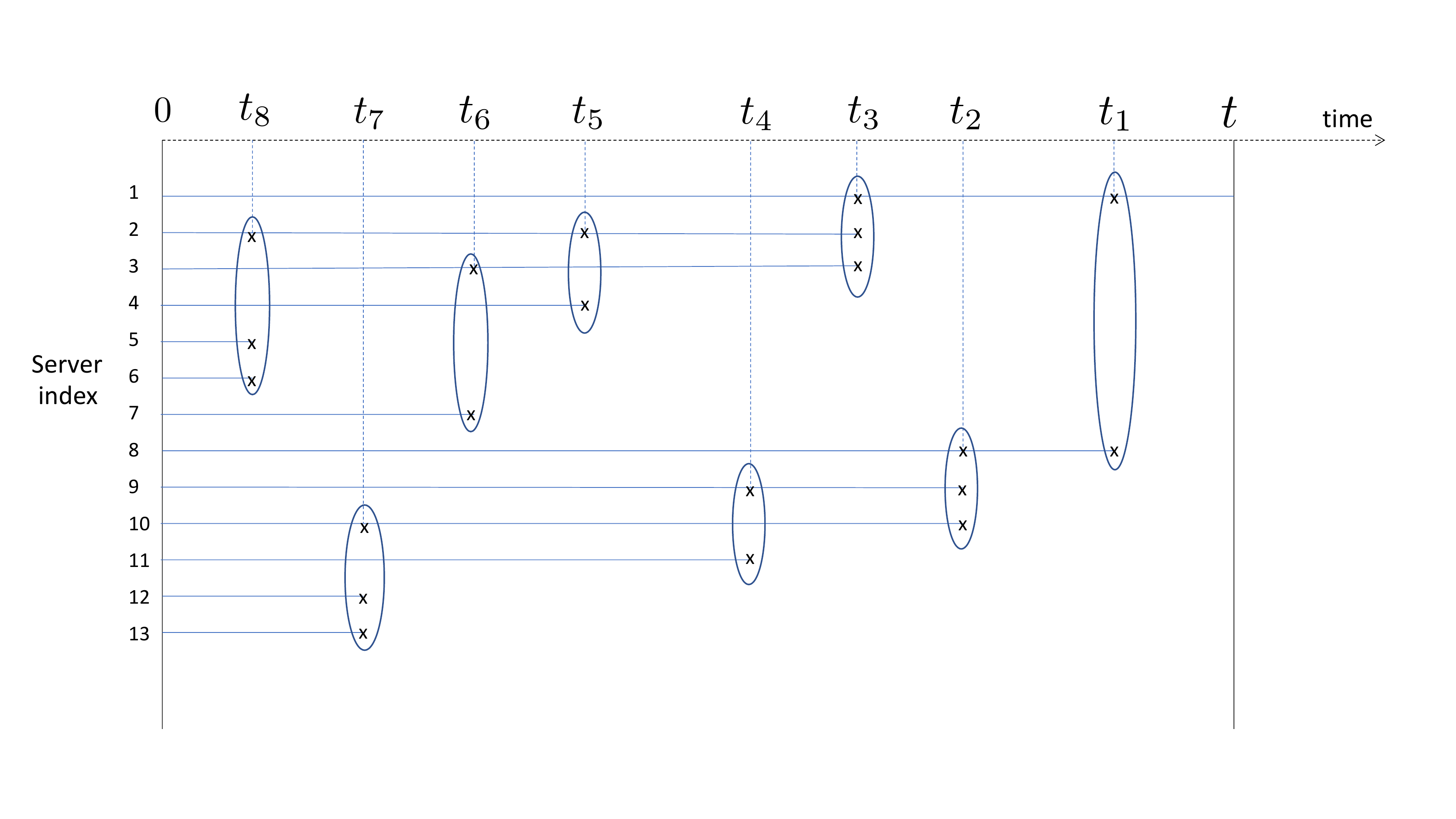}
\caption{Construction of dependence set $\bar D_1(t)$.}
\label{fig1}
\end{figure}

Fix $t\ge 0$. Suppose we consider a server, labeled $1$ to be specific. Let $U_1(t)$ denote its workload at time $t$.
Just like for our original system (with a finite number of servers), we will use the terminology of a job {\em selection set}, although here the latter will be defined formally, not as a result of an actual selection process.
Denote $\alpha_j = \lambda_j d_j$, $\alpha=\sum_j \alpha_j$. 
Then, by definition, the job arrivals of type $j$ selecting server 1 occur according to an independent Poisson process of rate $\alpha_j$.
We now define the dependence set $\bar D_1(t)$ of server $1$ at time $t$. It is defined via a branching process, running in reverse time 
from $t$ to $0$; the reverse time index is $0 \le s \le t$, with $s$ corresponding to actual (forward) time $t-s$.
To improve the exposition, we will define the construction via an example, shown in Figure~\ref{fig1};
in this example we also assume that there are two job classes, with $d_1=3$ and $d_2=2$.
Ovals indicate job arrivals, with crosses showing the servers they select.
By definition, set $A(s)$ at $s=0$ is $A(0)=\{1\}$. As $s$ increases, set $A(s)$ increases as follows. 
There is a Poisson process (in reversed time $s$) of rate $\alpha_j$, for each class $j$, of job arrivals selecting server $1$.
Suppose first such arrival occurs at time $s_1$ (which corresponds to forward time $t_1=t-s_1$ shown in  Figure~\ref{fig1}), and 
it happens to be a class-2 job arrival. The set $A(s)$ remains unchanged in the interval $[0, s_1]$, and at time $s_1$ we add $d_2-1=1$ servers to the set, which we call ``children'' of node $1$ added at time $s_1$; each added child server is a new server with distinct index -- it is server $8$ in Figure~\ref{fig1}. So, $A(s_1+) = \{1,8\}$, if we adopt a (non-essential) convention that $A(s)$ is left-continuous.
Starting time $s_1$, each server in $A(s)$ receives job arrivals according to independent Poisson processes, same as for server $1$.
Let $s_2$ (corresponding to forward time $t_1=t-s_1$ in  Figure~\ref{fig1}) be the next job arrival time, to any of the servers in set $A(s)$;
in Figure~\ref{fig1} it is class-1 job arrival to server $8$. Then $A(s)$ is constant in $(s_1,s_2]$, and at time $s_2$ we add $d_1-1=2$
new servers, with distinct indices $9,10$ -- these are children of server $8$ added at $s_2$. After $s_2$ we consider independent Poisson job arrival processes to each of servers in $A(s)$. The next job arrival is time $s_3$ (corresponding to forward time $t_3=t-s_3$ in  Figure~\ref{fig1}), which happens to be of class $1$ at server $1$, which adds two new children, $2$ and $3$, of server $1$ at time $s_3$;
and so on. In Figure~\ref{fig1} there are $8$ job arrivals in total in $[0,t]$, and the final set $A(t)= \{1,2,\ldots,13\}$.
We emphasize that {\em each added child server in this construction is a new server with a distinct index.} 
Then, by definition, $\bar D_1(t) = A(t)$. From now on, with a slight abuse of terminology, when we say ``dependence set  $\bar D_1(t)$,'' we assume that $\bar D_1(t)$ includes not only the set of servers, but also the graph describing the ``parent-child'' dependence structure and the job arrival times and classes. (For $\bar D_1(t)$, by construction, the graph describing the dependence structure is a tree.)

Now, given a realization of $\bar D_1(t)$, the random value of $U_1(t)$ is obtained by letting the initial workloads of all servers $i\in \bar D_1(t)$
to be i.i.d. with the distribution $\pr\{U_i(0) = w_k\} = a_k$, $k=1,\ldots,K$, and with the component size vectors for the involved  job arrivals
being independent with the corresponding distributions. As usual, between the times of job arrivals selecting a server, the workload $U_i(t)$ of each server decreases at rate $1$ (unless and until it reaches $0$). Clearly, if a server $i \in \bar D_1(t)$ is added as a child
at (forward) time $\tau$, the evolution of its workload after time $\tau$ does not affect the value of $U_1(t)$.
This completes the definition of $U_1(t)$.

\begin{lem}
\label{lem-graph-rep1}
For any $t\ge 0$, the (random) cardinality $f(t)=|\bar D_1(t)|$ of the dependence set $\bar D_1(t)$ is finite. 
Moreover, $f(t)$ satisfies
\beql{eq-dep-set-ode}
f'(t) = \gamma f(t), ~~\mbox{with}~ \gamma= \sum_j \alpha_j (d_j -1),
\eeql
and therefore
$$
\E |\bar D_1(t)| = e^{\gamma t}, ~t\ge 0.
$$
\end{lem}

 {\em Proof.} The proof uses the branching process argument. (It is analogous to that used in 
 %the proof of Lemmas 3.1 and 3.2 in \cite{GMP97} or in 
 Section 7 of \cite{BLP2012-jsq-asymp-indep}.) 
 Clearly, $f(s), ~s\ge 0,$ can be equivalently viewed as the cardinality of set $A(s)$ in the definition of $\bar D_1(t)$.
 In a small time interval $[s,s+\Delta s]$, the expected total number of all children
 that will be added is $f(s)[\sum_j \alpha_j (d_j -1) \Delta s + o(\Delta s)]$. This leads to ODE \eqn{eq-dep-set-ode}.
 We omit further details.
 $\Box$

\begin{lem}
\label{lem-graph-rep2}
For any $t\ge 0$, as $n\to\infty$,
\beql{eq-basic-mf-limit-single-queue-weak} 
W_1^n(t) \Rightarrow U_1(t)
\eeql
and, moreover,
\beql{eq-basic-mf-limit-single-queue} 
\pr\{ W_1^n(t) > w\}  \rightarrow \pr\{ U_1(t) >w\} ~~\mbox{for any $w\ge 0$}.
\eeql
\end{lem}

{\em Proof.}  Note that the workloads $W_i^n(t)$ are those of the servers $i=1,\ldots,n$, in the system with finite $n$.
(Also recall that initially the servers' indices are assigned in a random order.) 
Let us define the dependence set $D^n_i(t)$ for server $i$ at time $t$, in the system with given $n$. 
The construction of $D^n_i(t)$ is analogous to the construction of $\bar D_1(t)$ for the ``infinite system,'' except
here, as we increase the set $A^n_i(s)$ (analogous to set $A(s)$ in the definition of $\bar D_1(t)$) in reversed time $s$, the children are added according to actual job arrivals in the finite system with $n$ servers. As a result: (a) a single job arrival may simultaneously select multiple servers in $A^n_i(s)$; (b) when children are added at time $s$, they may be servers already present in $A^n_i(s)$; (c) the graph describing the parent-child dependence of servers in $D^n_i(t)$ is not necessarily a tree. 

Note that, marginally for each server $i$, 
the job arrival processes of different classes $j$ selecting this server are still independent Poisson processes with rates $\alpha_j$;
however, the independence of these arrival processes across the servers no longer holds. On the other hand, consider any fixed $C>0$.
Suppose there is a job arrival selecting at least one server in $A^n_i(s)$ at (reversed) time $s$.
Consider event 
$E(A^n_i(s)) =$   $\{$The arriving job
selects more than one server $A^n_i(s)$ or any of the children servers added by the arrival are already in $A^n_i(s)$.$\}$
Then, as $n\to\infty$, 
$\pr\{E(A^n_i(s))\} \to 0$, uniformly in all possible $A^n_i(s)$ with $|A^n_i(s)|\le C$. This in turn implies that, as $n\to\infty$,
 the rate at which a class-$j$ job arrival selecting a server in $A^n_i(s)$ occurs, converges to $|A^n_i(s)|\alpha_j$, 
 uniformly in $A^n_i(s)$ with $|A^n_i(s)|\le C$.
 
 Using the above observations, the proof proceeds by showing that, as $n\to\infty$, the set $D^n_i(t)$ converges to $\bar D_1(t)$, in an appropriate sense specified below.
%The formal proof proceeds as follows. 
It is easy to see that, for all sufficiently large $n$, 
$$
|D^n_1(t)| \le_{st} |\bar D_{1,\epsilon}(t)|,
$$
where $\bar D_{1,\epsilon}(t)$ is constructed the same way as $\bar D_{1}(t)$, but with arrival rates $\lambda_j$ replaced by slightly larger rates $(1+\epsilon)\lambda_j$, $\epsilon>0$.
By Lemma~\ref{lem-graph-rep1}, 
\beql{eq-dep-set-cardin}
%\E |D^n_1(t)| 
\E |\bar D_{1,\epsilon}(t)| = e^{\gamma' t}, ~t\ge 0,
\eeql
where $\gamma'= (1+\epsilon) \sum_j \alpha_j (d_j -1)$. We conclude that, uniformly in all sufficiently large $n$, 
$|D^n_1(t)|$ is stochastically dominated by a proper non-negative random variable. Using this fact, along with the observations in the previous paragraph,
we can easily couple the constructions of $D^n_1(t)$ for each $n$ and 
the construction of $\bar D_1(t)$ in such a way that,  w.p.1, 
\beql{eq-dep-sets-converge}
D^n_1(t) \to \bar D_1(t), ~~W_1^n(t) \to U_1(t),
\eeql
where the $D^n_1(t) \to \bar D_1(t)$ is defined as follows. We say that (the realizations of) $D^n_1(t)$ and $\bar D_1(t)$ are equivalent, and write $D^n_1(t) = \bar D_1(t)$ if they are equal, {\em up to relabeling of servers other than server $1$}; this means that they have 
equal sets of servers, equal number of job arrivals, equal parent-child dependence graph, and equal job arrival times.
We say that the convergence $D^n_1(t) \to \bar D_1(t)$ holds (for the realizations of $D^n_1(t)$ and $\bar D_1(t)$) if $D^n_1(t) = \bar D_1(t)$ for all sufficiently large $n$. The convergence \eqn{eq-dep-sets-converge} proves \eqn{eq-basic-mf-limit-single-queue-weak}.
Moreover, clearly the coupling construction can be such that, w.p.1, for all sufficiently large $n$, we in fact have $W_1^n(t) = U_1(t)$;
this proves \eqn{eq-basic-mf-limit-single-queue}.
$\Box$

Let us denote: 
\beql{eq-x-limit-def}
x^c_w(t) \doteq \pr\{U_1(t) > w\}.
\eeql

\begin{lem}
\label{lem-fsp-conv}
For the sequence of systems, considered in this section, the following holds for any 
$w\ge 0$ and $t\ge 0$:
\beql{eq-basic-mf-limit-gen} 
x_w^{n,c}(t) \stackrel{P}{\rightarrow} x_w^{c}(t).
\eeql
\end{lem}

{\em Proof.}  Consider a fixed $t$, a fixed $w$, and two fixed servers 1 and 2. Denote by $E^n_i \doteq \{W^n_1(t) > w\}$, $i=1,2$. To prove \eqn{eq-basic-mf-limit-gen}, it suffices
to prove that the covariance of $\bI(E^n_1)$ and $\bI(E^n_2)$ vanishes, namely
\beql{eq-cov-vanish}
\E \bI(E^n_1,E^n_2) - \E \bI(E^n_1) \E \bI(E^n_2) \to 0, ~~n\to\infty.
\eeql
Indeed, by the symmetry (exchangeability of the servers), from \eqn{eq-basic-mf-limit-single-queue} and \eqn{eq-x-limit-def} we have that $\E x_w^{n,c}(t) = \E \bI(E^n_1) = \pr\{W^n_1(t) > w\} \to \pr\{U_1(t) > w\} = x_w^{c}(t)$;
then the vanishing covariance \eqn{eq-cov-vanish} implies that the variance of $ x_w^{n,c}(t)$ vanishes as well.

In turn, to prove \eqn{eq-cov-vanish} it suffices to prove
\beql{eq-cov-vanish2}
\E \bI(E^n_1,E^n_2) \to  [\pr\{U_1(t) > w\}]^2, ~~n\to\infty.
\eeql
This follows from a straightforward extension of the constructions we used in the proof of Lemma~\ref{lem-graph-rep2}.
Namely, consider two independent constructions of the ``infinite server system'' dependence sets $\bar D_1(t)$ and $\bar D_2(t)$;
the latter is defined the same way as the former, except the ``seed'' server is $2$, not $1$; 
WLOG we can assume that the sets of server indices in $\bar D_1(t)$ and $\bar D_2(t)$ do not intersect.
Then consider the joint construction of the dependence sets $D^n_1(t)$ and $D^n_2(t)$, which starts with the ``seed''
set $A^n(0) = \{1,2\}$, corresponding to time $t$, and then goes ``backwards in time.'' In this case $D^n_1(t)$ and $D^n_2(t)$ are
{\em not} independent, but exactly the same argument as in the proof of Lemma~\ref{lem-graph-rep2} shows
that a coupling exists such that, w.p.1,  $(D^n_1(t),D^n_2(t)) \to (\bar D_1(t),\bar D_2(t))$
and $(W_1^n(t),W_2^n(t)) \to (U_1(t),U_2(t)$.
(The convergence $(D^n_1(t),D^n_2(t)) \to (\bar D_1(t),\bar D_2(t))$ is defined the same way as $D^n_1(t) \to \bar D_1(t)$, except the equivalence $(D^n_1(t),D^n_2(t)) = (\bar D_1(t),\bar D_2(t))$ is understood as equality up to relabeling of the servers other than $1$ and $2$.)
Recall $\bar D_1(t)$ and $\bar D_2(t)$ (and $U_1(t)$ and $U_2(t)$) are independent.
Moreover, the coupling can be such that, w.p.1, for all sufficiently large $n$, we in fact have 
$(W_1^n(t),W_2^n(t)) = (U_1(t),U_2(t)$. This proves \eqn{eq-cov-vanish2}, and then the lemma.
$\Box$

\section{Fluid sample paths}
\label{sec-fsp}

Suppose we are in the setting of Section~\ref{sec-asymp-indep-finite-int}.
Defined there function $x^{c}(t), t\ge 0,$ with values in $\cx$, we will call a {\em fluid sample path} (FSP). 
Clearly, an FSP initial state is: $x_w^{c}(0) = \sum_k a_k \bI\{a_k >w\}$, $0\le w < \infty$.
(Note that an FSP, by definition, arises as a result of the limiting procedure specified in Section~\ref{sec-asymp-indep-finite-int}.
Namely, the initial states of the pre-limit systems are such that exactly a fraction $a_k$ of servers has workload exactly $w_k$,
for some fixed parameters $a_k>0$ (such that $\sum_k a_k = 1$) and $w_k \in [0,\infty]$. In this paper we will only 
need FSPs defined this way.)

By $x^{c,\emptyset}(t), t\ge 0,$ we will denote the special FSP with $a_1=1, w_1=0$; this means that each pre-limit system starts from the ``empty'' initial state, with all initial workloads being $0$. Of course, $x_0^{c,\emptyset}(0) = 0$. This is the FSP ``starting from the empty initial state.'' As a special case of Lemma~\ref{lem-fsp-conv}, we obtain that for any fixed $c\le \infty$,  $w\ge 0$, $t\ge 0$,
\beql{eq-basic-mf-limit} 
x_w^{n,c,\emptyset}(t) \stackrel{P}{\rightarrow} x_w^{c,\emptyset}(t).
\eeql

The FSP definition and Lemma~\ref{lem-monotone1} imply the following monotonicity property for the FSPs.

\begin{lem}
\label{lem-fsp-monotone}
(i) Consider two FSPs, 
 $x^{c}(\cdot)$ and $\hat x^{\hat c}(\cdot)$,
such that $c \le \hat c$ and $x^{c}(0) \le \hat x^{\hat c}(0)$. Then $x^{c}(t) \le \hat x^{\hat c}(t)$ for all $t\ge 0$.
(ii) Consider two FSPs, 
 $x(\cdot)$ and $\hat x(\cdot)$, such that, for some $0\le h <\infty$ and $\tau \ge 0$, $\hat x_h(0)=0$ 
 and $x_h(\tau) > \hat x_0(0)$. Then, 
$\hat x(t) \le x(\tau+t)$ for all $t\ge 0$.
\end{lem}

\section{Properties of FSPs starting from the empty initial state}
\label{sec-fsp-from-empty}

In this section we study the properties of the FSPs $x_w^{c,\emptyset}(t), ~t\ge 0,$ starting from the empty initial state. 
Recall that $c\in [0,\infty]$ is the truncation parameter. Note that an FSP $x_w^{c,\emptyset}(t), ~t\ge 0,$ can be viewed as 
a scalar function of $(c,w,t)$.

\begin{lem}
\label{lem-limit-basic-monotone}
Function $x_w^{c,\emptyset}(t)$ is non-decreasing in $c,t$, and is non-increasing in $w$.
\end{lem}

 {\em Proof.}  This follows from \eqn{eq-basic-mf-limit}, along with Lemma~\ref{lem-monotone1} and Corollary~\ref{cor:zero}.
$\Box$

\begin{lem}
\label{lem-limit-basic-proper-decreasing}
For any $c>0$ and $t>0$,  $x^{c,\emptyset}(t)$ is proper and $x_w^{c,\emptyset}(t)>0$ for each $w<c$.
\end{lem}

 {\em Proof.}  By definition \eqn{eq-x-limit-def}, $x_w^{c,\emptyset}(t)=\P\{U_1(t) > w\}$, where the construction 
 of $U_1(t)$ is for the special case when all servers' workloads are initially $0$.
 It follows from the construction that $U_1(t)$ is finite w.p.1 for any finite $t\ge 0$, that is 
 $x^{c,\emptyset}(t)$ is proper. 
 %(If $c<\infty$ this is, of course, automatic.) 
 Moreover, it is easy to see from the construction that, for any $t\ge 0$ and $w<c$, the 
 probability $\P\{U_1(t) > w\}$ must be positive.
 $\Box$

 Sometimes, as in the proof of the next lemma, it will be convenient to interpret a given server workload evolution as the movement of a ``particle'' in $[0,\infty]$, with the workload being the particle location. With this interpretation, between the times 
of job arrivals that select the server, the particle moves left at the constant speed $1$ until/unless it ``hits'' $0$. 
At the times when a new job arrival adds to the server workload, the particle ``jumps right'' by the distance equal to the added workload.

\begin{lem}
\label{lem-limit-lipschitz}
As a function of $w\ge 0$, $x_w^{c,\emptyset}(t)$ is Lipschitz, uniformly in $c\le \infty$ and $t\ge 0$.
\end{lem}

 {\em Proof.}  Consider time $t$ and interval $[w,w+\delta]$. Consider the (pre-limit) system and process, with fixed $n$.
 All particles (server workloads) that are in $[w,w+\delta]$ at time $t$, at time $t+\delta$ will be in $[0,w]$, unless they are selected by new job arrivals in $[t,t+\delta]$. 
 Recall that new jobs arrive as a Poisson process of (unscaled) rate $\lambda n$, for a given $n$, and each job selects at most $d$ particles.
 Let $g^n_w$ be the (scaled) number of particles that cross point $w$ from left to right in the interval $[t,t+\delta]$.
  By the law of large numbers, 
 $$
 (P) \!\limsup_{n\to\infty} g^n_w \le \lambda d \delta.
 $$
 Therefore,
 $$
 x_w(t+\delta) - x_w(t) = (P) \!\lim_{n\to\infty} [x_w^n(t+\delta) - x_w^n(t)] \le 
 $$
 $$
 - (P) \!\lim_{n\to\infty} [x_w^n(t) - x_{w+\delta}^n(t)] + (P) \!\limsup_{n\to\infty} g^n_w \le 
 - [x_w(t) - x_{w+\delta}(t)] + \lambda d\delta .
 $$
 But, $0 \le x_w(t+\delta) - x_w(t)$.
 Therefore, $x_w(t) - x_{w+\delta}(t) \le \lambda d\delta$.
$\Box$

\begin{lem}
\label{lem-limit-incr-in-t}
For any fixed $0\le w < \infty$, function $x_w^{\emptyset}(t), t\ge 0,$ is 
strictly increasing in $t$. (Note that here we consider specifically the non-truncated system, $c=\infty$.)
\end{lem}

{\em Proof.} 
Fix any $\tau \ge 0$, any $\delta>0$, and any $w \ge 0$. We will show that $x^{\emptyset}_w(\tau+\delta) > x^{\emptyset}_w(\tau)$. Suppose, first, that $w>0$. By Lemma~\ref{lem-limit-basic-proper-decreasing},
at time $\delta$, $x^{\emptyset}(\delta)$ is such that $x^{\emptyset}_u(\delta)$ is positive for all $u<\infty$. 
Consider a state, let us denote it $\hat x(\delta)$, such that $\hat x_u(\delta) = a \in (0,1)$ for $u<w+\tau$, and 
$\hat x_u(\delta) = 0$ for $u\ge w+\tau$. (Fraction $a$ of  ``servers'' have workload exactly $w+t$, while the rest of the ``servers'' have workload $0$.) We can and do pick $a>0$ small enough so that $x^{\emptyset}_u(\delta) > \hat x_u(\delta)$ for all $u\ge 0$.
For each $n$, let us consider the process $x^{n,\emptyset}_u(\cdot)$ in the time interval $[\delta,\infty)$,
and compare it with the process $\hat x^n(\cdot)$ in the same time interval, starting from state $\hat x^n(\delta)=\hat x(\delta)$.
We have that
$$
\lim_{n\to\infty} \P\{x^{n,\emptyset}_u(\delta) \ge \hat x^{n}_u(\delta), ~\forall u \ge 0\} = 1.
$$
Given this, and using the monotonicity, we can couple these processes in a way such that
$$
\lim_{n\to\infty} \P\{x^{n,\emptyset}_u(\delta+t) \ge \hat x^{n}_u(\delta+t), ~\forall u \ge 0,  ~\forall t \ge 0\} = 1.
$$
By Lemma~\ref{lem-fsp-conv}, for any $t\ge 0$ and $u \ge 0$, 
$\hat x^{n}_u(\delta+t) \stackrel{P}{\rightarrow} \hat x_u(\delta+t)$,
where, by monotonicity, $\hat x_u(\delta+t) \ge x^{\emptyset}_u(t)$.
Now, by Lemma~\ref{lem-jump-remains}, 
$$
(P) \!\liminf_{n\to\infty} (\hat x^n_{w+\tau-t-}(\delta+t) - \hat x^n_{w+\tau-t}(\delta+t), ~0 \le t < w+\tau) \ge a e^{-\lambda d^2 (w+\tau)}=\epsilon >0.
$$
This implies that, for any $0 \le t < w+\tau$,
$$
\hat x_{w+\tau-t-}(\delta+t) - \hat x_{w+\tau-t}(\delta+t) \ge \epsilon,
$$
and then
$$
\hat x_{w+\tau-t-}(\delta+t) - x^{\emptyset}_{w+\tau-t}(t) \ge \epsilon,
$$
and finally
\beql{eq-jump-lower-bound}
x^{\emptyset}_{w+\tau-t}(\delta+t) - x^{\emptyset}_{w+\tau-t}(t) \ge \epsilon.
\eeql
In particular, substituting $t=\tau$, we obtain 
$x^{\emptyset}_{w}(\delta+\tau) - x^{\emptyset}_{w}(\tau) \ge \epsilon > 0$, which completes the proof for the case $w>0$.
To treat the case $w=0$, observe that the proof of \eqn{eq-jump-lower-bound} in fact holds as is, 
with the same $\epsilon>0$ in the RHS,
if in the LHS we replace $w$ by any $u \in (0,w]$. Therefore,
$x^{\emptyset}_{u}(\delta+\tau) - x^{\emptyset}_{u}(\tau) \ge \epsilon > 0$ for all small positive $u>0$, 
and then, by the continuity of $x^{\emptyset}_{u}(\tau)$, for $u=0$ as well.
$\Box$

As a corollary from the results in this section, we obtain the following uniform convergence to an FSP starting from the empty initial state.
\begin{lem}
\label{lem-conv-to-fsp}
For any $t\ge 0$,
$\| x^{n,c,\emptyset}(t) - x^{c,\emptyset}(t)\| \stackrel{P}{\rightarrow} 0$.
\end{lem}

\section{Fixed point}
\label{sec-fp}

Given the properties of the FSPs derived above, we see that,
as $t\to\infty$,
$$
x^{c,\emptyset}(t) \stackrel{u.o.c.}{\rightarrow} x^{*,c} \in \cx.
$$
The element $x^{*,c} \in \cx$ we will call the {\em fixed point} (for a given $c$). In particular, $x^*=x^{*,\infty}$. 

Given that all functions $x_w^{c,\emptyset}(t), ~w\ge 0,$ are uniformly Lipschitz, their u.o.c. limits $x_w^{*,c}, ~w\ge 0,$
are uniformly Lipschitz as well.

\begin{lem}
\label{lem-fp-load-upper}
The fixed point $x^{*,c}$ is such that $x_0^{*,c} \le \lambda$.
\end{lem}

{\em Proof.}  The proof is by contradiction. Suppose not, i.e. $x_0^{*,c} > \lambda$. Choose $T>0$ sufficiently large, so that 
$$
(1/T) \int_0^T x_0^{c,\emptyset}(t) dt > \lambda.
$$
Note that, for any $n$,
$$
(1/T) \int_0^T \E x_0^{n,c,\emptyset}(t) dt \le \lambda,
$$
because the LHS is the expected amount of work (per server per unit time) processed by the system in $[0,T]$ -- it cannot exceed $\lambda$,
which is the expected amount of work  (per server per unit time) arrived into the system in $[0,T]$.
But we have
$$
\lim_{n\to\infty} (1/T) \int_0^T \E x_0^{n,c,\emptyset}(t) dt = (1/T) \int_0^T x_0^{c,\emptyset}(t) dt > \lambda.
$$
This contradiction completes the proof.
$\Box$

\begin{lem}
\label{lem-up-lim-of-trunk}
For any $w \in [0,\infty)$, as $c \uparrow \infty$, $x^{*,c}_w \uparrow x^{*}_w$.
\end{lem}

 {\em Proof.} It follows from Lemma~\ref{lem-equality-persists}, that the processes with a finite $c$ and $c=\infty$ can be coupled so that
 $x_w^{n,c,\emptyset}(t) = x_w^{n,\emptyset}(t)$ for $t < c-w$. Therefore, $x_w^{c,\emptyset}(t) = x_w^{\emptyset}(t)$ for $t < c-w$.
 Then, for any $\epsilon>0$, we can choose a sufficiently large $t$ and then sufficiently large $c$, so that
 $$
 x_w^{c,\emptyset}(t) = x_w^{\emptyset}(t) > x_w^* - \epsilon,
 $$
 and then $x^{*,c}_w > x_w^* - \epsilon$.
 $\Box$

In particular, from Lemma~\ref{lem-up-lim-of-trunk}, $x^{*,c}_0 \uparrow x^{*}_0$ as $c \uparrow \infty$.

 Note that, for a finite $c$, the fixed point $x^{*,c}$ is automatically proper, since $x_c^{*,c}=0$.
 Our next goal is to show that $x^*$ is proper, i.e. $x_w^{*} \downarrow 0$ as $w\to\infty$, and $x_0^*=\lambda$.
 
 \begin{lem}
\label{lem-fp-proper-to-conserv}
If fixed point $x^*$ is proper, then $x^*_0 = \lambda$.
\end{lem}

 {\em Proof.}  Fix any $\epsilon > 0$ and any $\epsilon'>0$. Then fix a large $h>0$ such that
 $$
 1-(1-x_h^*)^d < \epsilon';
 $$
 if we have a subset consisting of $(1-x_h^*)n$ servers, then the probability that an arriving job selects a server outside this set is less than $\epsilon'$. Finally, fix a sufficiently large finite $c>h$ so that if a job selects only servers with workloads at most $h$, then the expected amount of (this job's) work lost due to truncation is less than $\epsilon$.
 
 For each $n$, consider process $x^{n,c,\emptyset}(\cdot)$, with the chosen above truncation parameter $c$.
 For any $t\ge 0$, since $x_h^{c,\emptyset}(t) \le x_h^{\emptyset}(t) < x_h^*$,
 $$
  \lim_{n\to\infty} \P\{1-(1-x_h^{n,c,\emptyset}(t))^d < \epsilon' \} = 1.
 $$
 Let $g^n_{lost}(T)$ denote the expected total (scaled) amount of workload lost due to truncation in the interval $[0,T]$.
 Recall that the job arrivals process is Poisson. Then,
 $$
 \limsup_{n\to\infty} g^n_{lost}(T) [\lambda \epsilon' + \lambda (1-\epsilon') \epsilon] T \le \lambda (\epsilon' + \epsilon) T.
 $$
 Let $g^n_{out}(T)$ denote the expected total (scaled) amount of workload processed (and left the system) in the interval $[0,T]$.
 $$
 \lim_{n\to\infty} g^n_{out}(T) = \lim_{n\to\infty} \int_0^T \E x_0^{n,c,\emptyset}(t) dt =
 \int_0^T x_0^{c,\emptyset}(t) dt.
 $$
 Recall that the process starts from the empty state, so by work conservation, for any $T>0$,
 $$
 g^n_{out}(T) \ge \lambda T - g^n_{lost}(T) - c.
 $$
 We obtain
 $$
 \frac{1}{T}  \int_0^T x_0^{c,\emptyset}(t) dt \ge \lambda (1-\epsilon' - \epsilon) -c/T.
 $$
 Letting $T\to\infty$, and recalling that $\lim_{t\to\infty} x_0^{c,\emptyset}(t) =  x_0^{*,c}$,
 $$
 x_0^{*,c} \ge \lambda (1-\epsilon' - \epsilon). 
 $$
 Recall that such $c$ can be chosen for arbitrarily small $\epsilon' >0$ and $\epsilon>0$.
 Then $x_0^{*}= \lim_c x_0^{*,c} \ge \lambda$. And by Lemma~\ref{lem-fp-load-upper}, $x_0^{*} \le \lambda$. 
  $\Box$
 
We now make the following observation. Up to this point in the 
paper (except in the statements of the main results, Theorem~\ref{thm-main} and Corollary~\ref{cor-main}), we never used the condition $\lambda<1$.
In particular, the definition of the fixed point $x^{*,c}$ does not depend on the condition $\lambda<1$ and neither do the
statements and proofs of Lemmas~\ref{lem-fp-load-upper} and \ref{lem-fp-proper-to-conserv}. Using this fact, we 
obtain the following corollary from Lemmas~\ref{lem-fp-load-upper} and \ref{lem-fp-proper-to-conserv}, which will be used later.

\begin{cor}
\label{cor-proper}
Consider the dependence of the fixed point $x^*$ on $\lambda$. If for a given $\lambda$ the
fixed point $x^*$ is proper, then $x^*_0 = \lambda$ (and then necessarily  $\lambda\le 1$).
\end{cor}

It is not difficult to strengthen Corollary~\ref{cor-proper}, to
show that if $x^*$ is proper, then necessarily $x^*_0 = \lambda < 1$. But we will not need this fact.

From now on, some of our results/proofs do require that $\lambda<1$.

\begin{lem}
\label{lem-fp-decr}
Function $x_w^{*,c}, w\ge 0,$ is strictly decreasing for $0 \le w < c$. (This does require that $\lambda<1$.)
\end{lem}

 {\em Proof.} The proof is by contradiction. Let $u\ge 0$ be a point such that $x_w^{*,c}$ is flat in an interval 
 to the right of $u$, , i.e. $x_u^{*,c} = x_v^{*,c}=a$ for some $u < v < c$. Note that $a \le x_0^{*,c} \le \lambda < 1$.
 Pick any $\epsilon>0$, and then $t>0$ large enough, so that both $x_u^{c,\emptyset}(t)$ and $x_v^{c,\emptyset}(t)$ are in
 $(a-\epsilon, a)$. Pick any $0< b < 1-\lambda$. 
 At time $t$ let us consider the non-intersecting sets: $\bar B=\bar B(n)$ is the set of empty servers, $A=A(n)$ is the set of servers with workload 
 greater than $v$. For their (scaled) cardinalities, we know that
 $$
 \lim_{n\to\infty} \P\{1-x_0^{n,c,\emptyset}(t) > b\} = 1,
 $$
 $$
 \lim_{n\to\infty} \P\{x_v^{n,c,\emptyset}(t) > a-\epsilon\} = 1.
 $$
 Note that every server in set $A$ will have a workload greater than $u$ at time $t+(u-v)$.
 If event $\bar B \ge b n$ does hold at time $t$, let us pick a fixed subset $B=B(n)$ of $bn$ servers that are empty at $t$, 
 and let us keep track of the servers in set $B$ over the time interval $[t,t+(v-u)]$. By Lemma~\ref{lem-monotone1}
 the distribution of workloads within this set stochastically dominates that of the following process: we consider only the servers in $B$ and we ``ignore'' any job arrival which selects at least one server outside $B$. Such a lower bounding process has the same structure 
as our original process, except it has a smaller number of servers, $bn$, and the job type arrival rates per server, $\lambda'_j$, are different.
Applying Lemma~\ref{lem-limit-basic-proper-decreasing} and \eqn{eq-basic-mf-limit} to the lower bounding process, we obtain the following property. Denote by $g^n$ the (scaled) number of 
servers in $B$, which at time $t+(v-u)$ have workloads greater than $u$. Then, there exists $\delta >0$ such that
$$
 \lim_{n\to\infty} \P\{g^n > \delta\} = 1.
$$
Combining these estimates, we obtain that
$$
 \lim_{n\to\infty} \P\{x_u^{n,c,\emptyset}(t+(v-u)) > a-\epsilon+\delta\} = 1.
$$
But, this is true for any $\epsilon>0$. So, we must have
$
x_u^{c,\emptyset}(t+(v-u)) > a,
$
and then 
 $
x_u^{*,c} > a,
$
a contradiction.
 $\Box$

\begin{lem}
\label{lem-fp-proper}
Fixed point $x^*$ is proper and $x^*_0 = \lambda$. (This does require that $\lambda<1$.)
\end{lem}

 {\em Proof.}  Given Lemma~\ref{lem-fp-proper-to-conserv} (or Corollary~\ref{cor-proper}), it suffices to prove that $x^*$ is proper.
Suppose not, i.e. $x_w^* \downarrow = a > 0$ as $w\to\infty$. Then we claim the following. Consider the FSP $\hat x(t)$ corresponding to the 
initial state with fraction $a$ of servers having initially infinite workload and the remaining fraction $b=1-a$ of servers being empty.
(This is the system described in Section~\ref{sec-equiv}, with the B-subsystem being initially empty.) 
So, $\hat x_w(0)=a$ for all $w\ge 0$. We claim that
\beql{eq-conv-from-infinite-fraction}
\hat x_w(t) \uparrow x^*_w, ~\forall w \ge 0.
\eeql
(Consequently, $\|\hat x(t) - x^*\| \to 0$.) Note that by Lemma~\ref{lem-fsp-monotone}, 
$\hat x(t) \ge x^{\emptyset}(t)$, so $\lim_t \hat x_w(t) \ge x^*_w$. Then,
to prove \eqn{eq-conv-from-infinite-fraction}, it suffices to show that for any $t\ge 0$, $\hat x_w(t) \le x^*_w$.
For this, we consider the family of FSPs, $\hat x^{(h)}(\cdot)$, parameterized by $h< \infty$, corresponding to initial states such that the fraction $a$ of servers have initial workloads $h$ and the remaining fraction $b=1-a$ of servers are empty. For any fixed $h$, there exists a sufficiently 
large $\tau>0$ such that $x^{\emptyset}(\tau) > \hat x^{(h)}(0)$. (Here we use the fact that $x^*_w > a$ for all finite $w\ge 0$, because
$x^*_w$ is strictly decreasing in $w$.) Then, by Lemma~\ref{lem-fsp-monotone}, $\hat x^{(h)}(t) \le x^{\emptyset}(\tau+t) < x^*$ for all $t\ge 0$.
It remains to notice that for any fixed $w\ge 0$ and $t\ge 0$, if we choose $h > w+t$, then $\hat x_w(t) = \hat x_w^{(h)}(t)$.
Thus, claim \eqn{eq-conv-from-infinite-fraction} is proved.

Denote by $\bar x^*$ the element of $\cx$, defined by $\bar x_w^* = (x^*_w-a)/b, ~w\ge 0.$ By this definition, $\bar x^*$ is proper.
Claim \eqn{eq-conv-from-infinite-fraction} proves
 that $\bar x^*$ is nothing else but the fixed point for the $B$-subsystem, as defined in Section~\ref{sec-equiv},
starting from the empty state.
Since $\bar x^*$ is proper, by Corollary~\ref{cor-proper}, $\bar x_0^* = \rho_B$. By \eqn{eq-rho_B}, $\rho_B \ge \rho=\lambda$.
But then $x^*_0 = a + b \bar x_0^* \ge a + b \lambda > \lambda$. This contradicts Lemma~\ref{lem-fp-load-upper}.
 $\Box$

 \section{Proof of Theorem~\ref{thm-main}}
 \label{sec-proof-main}
 
 Consider the sequence $x^n(\infty)$. The sequence of their distributions is tight because $\cx$ is compact. 
 Consider any fixed subsequence along which $x^n(\infty) \Rightarrow x^{**}$, where $x^{**}$ is a random element in $\cx$.
 It will suffice to show that, $x^{**}=x^*$.
 
 \begin{lem}
\label{lem-distr-lower-bound}
$x^* \le_{st} x^{**}$.
\end{lem}

{\em Proof.} We can construct a stationary version of the process, $x^n(\cdot)$, and the process $x^{n,\emptyset} (\cdot)$ on 
a common probability space, so that $x_w^{n,\emptyset} (t) \le x_w^n(t)$ for all $w$ and $t$. We conclude that 
$x^{n,\emptyset}(t) \le_{st} x^n(t)$ for any $t\ge 0$. Recall that $\| x^{n,\emptyset}(t) - x^{\emptyset}(t)\| \stackrel{P}{\rightarrow} 0$ as $n \to \infty$.
This implies that, for any fixed $h>0$ and $t\ge 0$,
$$
(P) \!\liminf_{n\to\infty} (x_w^n(\infty), ~w \in [0,h]) \ge (x_w^{\emptyset}(t), ~w \in [0,h]).
$$
Since this is true for any $t$, and $x_w^{\emptyset}(t)$ is strictly increasing in $t$, we conclude that, 
for any fixed $h>0$ and $t\ge 0$ the stronger property holds: 
\beql{eq-load-lower-bound-interval}
\lim_{n\to\infty} \pr\{x_w^n(\infty) \ge x_w^{\emptyset}(t), ~\forall w \in [0,h)\} = 1.
\eeql

Recalling that $x_w^{\emptyset}(t)$ is continuous in $w$, observe that the subset 
$\{y\in \cx~|~ y_w \ge x_w^{\emptyset}(t), ~\forall w \in [0,h)\}$ is closed. Therefore,
$$
\pr\{x_w^{**} \ge x_w^{\emptyset}(t), ~\forall w \in [0,h)\} \ge \limsup_n 
\pr\{x_w^n(\infty) \ge x_w^{\emptyset}(t), ~\forall w \in [0,h)\} = 	1.
$$
Thus, for any fixed $h>0$ and $t\ge 0$,
$$
\pr\{x_w^{**} \ge x_w^{\emptyset}(t), ~\forall w \in [0,h)\} = 	1.
$$
It remains to recall that $x_w^{\emptyset}(t) \uparrow x_w^*$ as $t\to \infty$, to finally conclude that
$$
\pr\{x_w^{**} \ge x_w^*, ~\forall w \ge 0 \} = 	1.
$$
$\Box$

For future reference, note that Lemma~\ref{lem-distr-lower-bound} (or 
\eqn{eq-load-lower-bound-interval}) implies, in particular, 
\beql{eq-load-lower-bound}
(P) \!\liminf_{n\to\infty} x_0^n(\infty) \ge x_0^*.
\eeql
% COM because $\{x_0 \ge x_0^* - \epsilon\}$ is an open subset of $\cx$.

 \begin{lem}
\label{lem-distr-equality}
$x^{**} = x^*$.
\end{lem}

{\em Proof.} 
Suppose not, i.e. there exists $w>0$ and $a>0$ such that
$$
\pr\{x^{**}_w > x_w^* + a \} = \delta > 0.
$$
We will show that this leads to a contradiction. 
Let $x^n(\cdot)$ denote
a stationary version of the process.
The subset
$
\{y\in \cx~|~ y_w > x_w^* + a\}
$
is open, so
$$
\liminf_n \pr\{x^{n}_w(\infty) > x_w^* + a \} \ge \delta.
$$
Pick a sufficiently small $\epsilon>0$, such that $\epsilon < a e^{-\lambda d^2 w}$.
Pick a sufficiently large $\tau>0$ and then a sufficiently small $u>0$, such that $u\le w$ and $x_u^{\emptyset}(\tau) + \epsilon > x_0^* + \epsilon/2$.
Let $\hat x^{n} (t) \doteq x^{n,\emptyset} (\tau+t)$, so that $\hat x^{n} (0)$ is equal in distribution to $x^{n,\emptyset} (\tau)$.
Since $x_v^\emptyset(\tau) < x^*_v$ for all $v \le w$, 
we know from the argument in the proof of Lemma~\ref{lem-distr-lower-bound}, that,
if we take independent initial distributions of $x^n(\cdot)$ and $\hat x^n(\cdot)$ then
$$
\lim_n \pr\{ x^n_v(0) \ge \hat x^n_v(0), ~\forall v \le w\} = 1.
$$
Let us couple $x^n(\cdot)$ and $\hat x^n(\cdot)$ in the natural way. By Lemma~\ref{lem-equality-persists-gen}, if condition $\{ x^n_v(0) \ge \hat x_v^n(0), ~\forall v \le w\} $ does
hold, then condition $\{ x^n_v(t) \ge \hat x^n_v(t), ~\forall v \le w-t, ~\forall t \le w\}$ holds as well. 
Furthermore, coupled with these processes, let us consider the following further modification $\tilde x^n(\cdot)$ of the process 
$x^n(\cdot)$. The initial state $x^n(0)$ is replaced by the initial state $\tilde x^n(0)$, where
$$
\tilde x^n_u(0) = x^n_u(0), ~u < w,
$$
$$
\tilde x^n_u(0) = \hat x^n_u(0), ~u \ge w.
$$
Note that $\tilde x_{w-}^n(0)-\tilde x_{w}^n(0) \ge  x_{w}^n(0)-\tilde x_{w}^n(0)$. Therefore, $\tilde x^n(0)$ is such that 
$$
\liminf_n \pr\{ \mbox{At least (scaled) number $a$ of servers have workload exactly $w$}\} \ge \delta.
$$
Consider the three coupled processes over the time interval $t\in [0,w-u]$. We can make the following conclusions:
$$
(P) \!\liminf_{n\to\infty} \hat x_{u}^n(w-u) \ge x_u^{\emptyset}(\tau+w-u) > x_u^{\emptyset}(\tau)
$$
which implies
$$
(P) \!\liminf_{n\to\infty} \tilde x_{u}^n(w-u) > x_u^{\emptyset}(\tau).
$$
We also can use Lemma~\ref{eq-jump-remains}  to conclude that 
$$
\liminf_n \pr\{ \tilde x^n_{u-}(w-u) -  \tilde x^n_{u}(w-u) \ge \epsilon \} \ge  \delta.
$$
Recalling our choice of $\epsilon$ and $u$, we obtain that
$$
\liminf_n \pr\{ \tilde x^n_{u-}(w-u) \ge x^*_0 + \epsilon/2 \} \ge  \delta.
$$
We also have
$$
\lim_n \pr\{ x^n_{u-}(w-u) \ge \tilde x^n_{u-}(w-u) \} = 1.
$$
The last two displays imply that
$$
\liminf_n \pr\{ x^n_{0}(w-u) \ge x^*_0 + \epsilon/2\} \ge  \delta.
$$
Recall that $x^n(\cdot)$ is the stationary version of the process. Then,
$$
\liminf_n \pr\{ x^n_{0}(\infty) \ge x^*_0 + \epsilon/2 \} \ge  \delta.
$$
From this display and \eqn{eq-load-lower-bound},
$$
\liminf_n \E x^n_{0}(\infty)  \ge x^*_0 + \delta \epsilon/2 > x^*_0 = \lambda,
$$
which contradicts the conservation law $\E x^n_{0}(\infty) = \lambda$.
$\Box$

Lemmas~\ref{lem-distr-lower-bound} and \ref{lem-distr-equality} imply that $x^n(\infty) \Rightarrow x^{*}$, thus proving Theorem~\ref{thm-main}.

%\iffalse
%\bibliographystyle{acmtrans-ims}
%\bibliographystyle{apt}
\bibliographystyle{abbrv}
\bibliography{biblio-stolyar}

\begin{thebibliography}{10}

\bibitem{adan2018fcfs}
I.~Adan, I.~Kleiner, R.~Righter, and G.~Weiss.
\newblock Fcfs parallel service systems and matching models.
\newblock {\em Performance Evaluation}, 127:253--272, 2018.

\bibitem{ayesta2018unifying}
U.~Ayesta, T.~Bodas, and I.~M. Verloop.
\newblock On a unifying product form framework for redundancy models.
\newblock {\em Performance Evaluation}, 127:93--119, 2018.

\bibitem{ayesta2019redundancy}
U.~Ayesta, T.~Bodas, and I.~M. Verloop.
\newblock On redundancy-d with cancel-on-start aka join-shortest-work (d).
\newblock {\em ACM SIGMETRICS Performance Evaluation Review}, 46(2):24--26,
  2018.

\bibitem{Bramson2011-jsq-stabil}
M.~Bramson.
\newblock Stability of join the shortest queue networks.
\newblock {\em The Annals of Applied Probability}, 21:1568--1625, 2011.

\bibitem{BLP2012-jsq-asymp-indep}
M.~Bramson, Y.~Lu, and B.~Prabhakar.
\newblock Asymptotic independence of queues under randomized load balancing.
\newblock {\em Queueing Systems}, 71:247--292, 2012.

\bibitem{Foss-Chernova-1998}
S.~Foss and N.~Chernova.
\newblock On the stability of a partially accessible multi-station queue with
  state-dependent routing.
\newblock {\em Queueing Systems}, 29:55--73, 1998.

\bibitem{Harcol-Balter-2017}
K.~Gardner, M.~Harchol-Balter, A.~Scheller-Wolf, M.~Velednitsky, and
  S.~Zbarsky.
\newblock Redundancy-d: The power of d choices for redundancy.
\newblock {\em Operations Research}, 65(4):1078--1094, 2017.

\bibitem{gardner2015reducing}
K.~Gardner, S.~Zbarsky, S.~Doroudi, M.~Harchol-Balter, and E.~Hyytia.
\newblock Reducing latency via redundant requests: Exact analysis.
\newblock {\em ACM SIGMETRICS Performance Evaluation Review}, 43(1):347--360,
  2015.

\bibitem{GMP97}
A.~Greenberg, V.~Malyshev, and S.~Popov.
\newblock Stochastic model of massively parallel computation.
\newblock {\em Markov Processes and Related Fields}, 2:473 -- 490, 1997.

\bibitem{hellemans2019performance}
T.~Hellemans, T.~Bodas, and B.~Van~Houdt.
\newblock Performance analysis of workload dependent load balancing policies.
\newblock {\em Proceedings of the ACM on Measurement and Analysis of Computing
  Systems}, 3(2):35, 2019.

\bibitem{PTW}
G.~Pang, R.~Talreja, and W.~Whitt.
\newblock Martingale proofs of many-server heavy-traffic limits for markovian
  queues.
\newblock {\em Probability Surveys}, 4:193--267, 2007.

\bibitem{shah2015redundant}
N.~B. Shah, K.~Lee, and K.~Ramchandran.
\newblock When do redundant requests reduce latency?
\newblock {\em IEEE Transactions on Communications}, 64(2):715--722, 2015.

\bibitem{St2014_pull}
A.~L. Stolyar.
\newblock Pull-based load distribution in large-scale heterogeneous service
  systems.
\newblock {\em Queueing Systems}, 80(4):341--361, 2015.

\bibitem{St2015_pull}
A.~L. Stolyar.
\newblock Pull-based load distribution among heterogeneous parallel servers:
  the case of multiple routers.
\newblock {\em Queueing Systems}, 85(1-2):31--65, 2017.

\bibitem{Vulimiri-2013}
A.~Vulimiri, P.~Godfrey, R.~Mittal, J.~Sherry, S.~Ratnasamy, and S.~Shenker.
\newblock Low latency via redundancy.
\newblock CoNEXT 2013 - Proceedings of the 2013 ACM International Conference on
  Emerging Networking Experiments and Technologies, pages 283--294. Association
  for Computing Machinery, Jan. 2013.

\bibitem{VDK96}
N.~Vvedenskaya, R.~Dobrushin, and F.~Karpelevich.
\newblock Queueing system with selection of the shortest of two queues: an
  asymptotic approach.
\newblock {\em Problems of Information Transmission}, 32(1):20--34, 1996.

\end{thebibliography}


\begin{thebibliography}{1}

\bibitem{St2014_pull}
A.~L.~Stolyar.
\newblock Pull-based load distribution in large-scale heterogeneous service systems.
%\newblock Bell Labs Technical Memo, 2014. Submitted.
\newblock \emph{Queueing Systems}, 
 2015. 
%Vol.47, No.1, to appear.
%pp.381-408.
%\\
http://arxiv.org/abs/1407.6343

\end{thebibliography}
%\bibliography{bibliography}
%\fi

%\input{repl-fp.bbl}

\iffalse

\fi

\end{document}